\documentclass[12pt]{article}
\usepackage{mathptmx}
\usepackage{graphicx}
\usepackage{bm}
\usepackage{amssymb}
\makeatletter
\author{Robert P. C. de Marrais \footnote{Email address:  rdemarrais@alum.mit.edu} \\ \noindent
\emph{Thothic Technology Partners, P.O.Box 3083, Plymouth MA 02361}}

\title{Placeholder Substructures II:~~  Meta-Fractals, Made of Box-Kites, Fill Infinite-Dimensional Skies}

\begin{document}
\maketitle
\makeatother

\begin{abstract}
Zero-divisors (ZDs) derived by Cayley-Dickson Process (CDP) from
N-dimensional hypercomplex numbers ($N$ a power of $2$, and at least
$4$) can represent singularities and, as $N \rightarrow \infty$,
fractals -- and thereby, scale-free networks.  Any integer $> 8$ and
not a power of $2$ generates a meta-fractal or \textit{Sky} when it
is interpreted as the \textit{strut constant} (S) of an ensemble of
octahedral vertex figures called \textit{Box-Kites} (the fundamental
ZD building blocks). Remarkably simple bit-manipulation rules or
\textit{recipes} provide tools for transforming one fractal genus
into others within the context of Wolfram's Class 4 complexity.
\end{abstract}

\section{Introduction By Way of Reprise:  From Box-Kites to ETs}
The creation of $2^{N}$-dimensional analogues of Complex Numbers
(and it was not a trivial insight of 19th Century algebra that
legitimate analogs \textit{always} have dimension a power of $2$) is
handled by a now well-known algorithm called the Cayley-Dickson
Process (CDP). Its name suggests a compressed account of its
history:  for Arthur Cayley -- simultaneously with, but
independently of, John Graves -- jumped on Hamilton's initial
generalization of the 2-D Imaginaries to the 4-D Quaternions within
weeks of its announcement, producing -- by the method later
streamlined into Leonard Dickson's close-to-modern ``cookie-cutter''
procedure -- the 8-D Octonions. The hope, voiced by no less than
Gauss, had been that an infinity of new forms of Number were lurking
out there, with wondrous properties just awaiting discovery, whose
magical utility would more than compensate for the loss of things
long taken for granted as their seekers ascended into higher
dimensions.  But such fantasies were quashed quite abruptly by
Adolph Hurwitz's proof, just a few years before the 20th Century
loomed, that it only took four dimension-doublings past the Real
Number Line to find trouble: the 16-D Sedenions had zero-divisors,
which meant division algebra itself broke down, which meant
researchers were so at a loss to find anything good to say about
such Numbers that nobody bothered to even give their 32-D immediate
successors a name, much less investigate them seriously.

But it is with these 32-D ``Pathions'' (for short for
``pathological,'' which we'll call them from now on) that our own
account will pick up in this second part of our study of
``placeholder substructures'' (i.e., ``zero divisors'') For, due to
a phenomenon we dubbed \textit{carrybit overflow} in the first
installment, strange yet predictable things are found to be afoot in
the ZD equivalent of a ``Cayley Table.''  As we'll see shortly, this
is a listing, in a square array, of the ZD ``emanations'' (or lack
of same) of all ZD ``elements'' with each other -- all, that is,
sharing membership in an ensemble defined not by a shared ``identity
element,'' but a common \textit{strut constant}.

What we'll see is that the lacks are of the essence: for each
doubling of $N$, the Emanation Table (ET) for the $2^{N+1}$-ions of
same strut-constant will contain that of its predecessor, leading to
an infinite ``boxes-within-boxes'' deployment whose empty cells
define, as $N$ grows ever larger, an unmistakable \textit{fractal
limit}. The full algorithmic analysis of such Matrioshka-doll-like
``meta-fractal'' aspects -- by the simple rules of what we'll call
``recipe theory'' (after the R, C, and P values related to the Row
label, Column label, and their cell-specific Products in such
Tables) -- must await our third and last installment.  But the
colored-quilt-like graphics can be viewed by any interested readers
at their leisure, in the Powerpoint slide-show online at Wolfram
Science from our mid-June presentation at NKS 2006.[1]  (The
slide-show's title is almost identical to that of this monograph, as
this latter is meant to be the ``theorem/proof'' exposition of that
iconic, hence largely intuitive and empirically driven narration.)

What we'll need to undertake this voyage is a quick reprise of the
results from Part I [2].  As the hardest part (as a hundred years of
denial would imply) is finding the right way to think about the
phenomenology of zero-division, not understanding its basic workings
once they're hit upon, such a summary can be much more brief and
easy to follow than the proofs required to produce and justify it.
We need but grasp 3 rather simple things.  First, we must
internalize the path and vertex structure of an Octahedron -- for,
properly annotated and storyboarded, this will provide us with the
Box-Kite representation that completely catalogs ZDs in the 16-D
arena where they first emerged (and, as we'll see in our Roundabout
Theorem herein, underwrites all higher-dimensional ZD emergences as
well).

Second, instead of the cumbersome apparatus of CDP that one finds in
algebra texts and the occasional software treatment, we offer two
easy algebraic one-liners which (inspired by Dr. Seuss's ``Thing 1''
and ``Thing 2''), we simply call ``Rule 1'' and ``Rule 2'' -- which
operate, in almost Pythagorean earnest, on triplets of integers
(indices of associative triplets among our Hypercomplex Units, as
we'll learn), and which, by so doing, accomplish everything the
usual CDP tactics do, but without the all-too-frequent obfuscation.
(There is also a very useful, albeit quite trivial, ``Rule 0,''
which merely states that any integer-triple serving to index an
associative triplet for one power of $N$ will continue to do so for
all higher powers.  What makes this useful is its allowing us to
recursively take triplet ``givens'' for lower-level $2^{N}$-ions
than those of current interest and toss them into the central circle
of the third thing we must grasp.)

We'll need, that is, to be able to draw the simplest finite
projective group's 7-line, 7-node representation, the so-called
PSL(2,7) triangle.  The Rules, plus the Triangle, applied to
Box-Kite edge-tracings and nodal indices, are all we'll need.
Indeed, the Box-Kite itself can be readily derived from the
Triangle, by suppressing the central node, and then recognizing four
correspondences.  First, see the Triangle's 3 triple-noded sides --
two vertices plus midpoint -- as the sources of the Box-Kite's trio
of ``filled-in'' triangles dubbed Trefoil Sails.
 Second, link the 1 triple-noded \textit{circle} (which \textit{is} a projective
\textit{line}, after all), wrapped around the suppressed center and
threading the midpoints, as the 4th such triangle, the quite special
``Zigzag Sail.''  Third, envision the 3 lines from midpoints to
angles as underwriting the ZD-challenged part of the diagram
(because ZDs housed at the midpoint node cannot mutually zero-divide
any housed at the opposite, vertex, node), the \textit{struts}
(whence strut \textit{constants}).  Fourth and last, imagine the
other four triangles of the Box-Kite (meeting, as with the first
four, each to each, at corners only, like same-colored checkerboard
squares) as the vents where the wind blows.  They keep the kite
afloat, letting the four prettily colored jib-shaped Sails show off,
while the trio of wooden or plastic dowels that form the struts
thanklessly provide the structural stability that makes the kite
able to fly in the first place.

As Euclid knew well, 3 points determine a Triangle as well as a
Circle -- which is how we can glibly switch gears between
representations based on these projective lines.  But the easy
convertibility of lines to circles is what \textit{projective} means
here -- and is, as well, at the very heart of linking the above
geometrical images to Imaginary Numbers.  From Argand's diagram to
Riemann's Sphere, this has been the essence of Complex geometry. On
the latter image only, place a sphere on a flat tabletop, call the
point of contact ${\bf S}$ (for ``South''), and then direct rays
from its polar opposite point ${\bf N}$.  Rays through the equator
intersect the table in a circle whose radius we ascribe an absolute
value of $1$, with center ${\bf S} = 0$.  This circle is just the
trace of the usual $e^{i \cdot 2 \pi \cdot \theta}$
exponential-orbit equation, with the $i$ in the exponent, of course,
being the standard Imaginary. Any diameter through this circle,
extended indefinitely in either direction, is clearly a ``projective
pencil'' of a circular motion in the plane containing both it and
${\bf N}$, and centered on the latter.

What each ``line,'' then, in the PSL(2,7) triangle represents is a
coherent system interrelating 3 distinct imaginaries, one per nodal
point:  that is, a ``Quaternion copy'' \textit{sans} the Reals
(which latter, like our ${\bf N, S}$ polar axis in the above, must
stand ``outside'' the Number Space itself, since 3-D visualization
is all used up by the nodes' dimensional requirements). Hence, the 7
lines are the 7 interconnected Quaternion copies which constitute
the 8-D Octonions. And what makes this especially rich for our
purposes is the built-in recursiveness of this Octonion-labeling
scheme for higher-dimensional isomorphs, embedded in the sorts of
ensembles we'll be needing ETs to investigate more thoroughly.

To see how this relates to actual integers, take the prototype of
the 7 lines in the Triangle, and consider the Quaternions strictly
from the vantage of CDP's Rule 1.  The first task in studying any
system of $2^{N}$-ions is generating its units, so start with $N =
0$. Treat this singleton as the index of the Real axis:  $i_{0}$,
that is, is identically $1$. Add a unit whose index = $2^{0}=1$ and
we have the complex plane.  Now, add in a unit whose index is the
next available power of $2$ -- with $N = 1$, this is $2$ itself.
Call this unit and its index ${\bf G}$ for \textit{Generator}, and
declare this inductive rule: the index of the product of any two
units is always the XOR of the indices of the units being
multiplied; but, for any unit with index $ u < {\bf G}$, the product
of said unit, written on the left (right), with the Generator
written on its right (left), has index equal to their indices'
simple sum, and sign equal (opposite) to the product of the signs of
their units':  $i_{1} \cdot i_{2} = + i_{3}$, but $i_{2} \cdot i_{1}
= - i_{3}$.  But this is just a standard way of summarizing
Quaternion multiplication.

Now, set N = 2, making ${\bf G} = 4$.  Applying the same logic, but
slightly generalized, we get three more triplets of indices.
Dispensing with the tedious overhead of explicitly writing the
indices as subscripts to explicit copies of the letter $i$, these
are written in cyclical positive order (CPO) as follows: $(1, 4, 5);
(2, 4, 6); (3, 4, 7)$.  (CPO is not mysterious:  it just means read
the triplet listing in left-right order, and so long as we multiply
any unit with any such index by the unit whose index is to the right
of it, the third term will result with signing as specified above:
e.g., $i_{4} \cdot i_{5} = + i_{1}$; $i_{4} \cdot i_{3} = - i_{7}$.)
We now have 4 of the Octonions' 7 triplets, forming labels on the
nodes of 4 of PSL(2,7)'s lines. Call the central circle spanning the
medians the Rule 0 line (the Quaternions' ``starter kit'' we just
fed into our Rule 1 induction machine). Putting ${\bf G}$ = 4 in the
center, the 3 lines through it are our Rule 1 triplets. If we
further array the Quaternion index-set $(1, 2, 3)$ in clockwise
order around the $4$, starting from the left slope's midpoint at 10
o'clock, these lines are all oriented pointing into the angles. Now,
with ``Rule 2,'' let's construct the lines along the Triangle's
sides.

Here's all that Rule 2 says:  given an associative index-triplet
(henceforth, \textit{trip}) like the Quaternions' $(1, 2, 3)$, fix
any one among them, then take its two CPO successors and add ${\bf
G}$ to them.  Swap the order of the resulting two new units, and you
have a new trip.  Hence, fixing 1, 2, and 3 in turn, in that order,
Rule 2 gives us these 3 triplets:  $( 1, 7, 6); (2, 5, 7); (3, 6,
5)$.  If you've drawn PSL(2,7) with the Octonion labels per the
instructions in the last paragraph, you've already seen these 3
trips are the answers ... and now you know how and why they're
oriented, too. (Clockwise, in parallel with the Rule 0 circle).

We've now laid out all the ingredients we need to do a basic
run-through of Box-Kite properties.  We'll merely state and describe
them, rather than prove them (but we'll give the Roman numerals of
the theorem numbers from last installment, for those who want to
follow them).  The first feature in need of elucidating, which
should have those who've been reading attentively scratching their
heads just about now, is this:  the relations between the indices at
the nodes of PSL(2,7) \textit{qua} Octonion labeling scheme are
clear enough; but how can these same labeled nodes serve to
underwrite the 16-D Sedenion framework that Box-Kites reside in? The
answer has two parts.

First part:  since all Imaginaries have negative Reals as squares,
Imaginaries whose products are zero must have different indices --
meaning that the simple case (which we call ``primitive'' ZDs) will
always involve products of \textit{pairs} of differently-indexed
units, whose respective planes share no points other than 0 [IV].
Second part: given any such ZD dyad, neither index can ever equal
${\bf G}$ [II]; and, one must have index $> {\bf G}$, while the
other has index $< {\bf G}$ [I, III]. The Octonion labeling scheme
maps to the four Sails of a legitimate Sedenion Box-Kite [V],
because it only provides the \textit{low}-index labels at each of
the 6 Octahedral vertices.

The $4$ in the center of our example, meanwhile, is no longer the
${\bf G}$ for this setup, since that role in now played by $8$ (the
next power of 2 in the CDP induction). In the context of the
Box-Kite scheme, it is now represented by a different letter: ${\bf
S}$, for \textit{strut constant} -- the only Octonion index
\textit{not} on a Box-Kite vertex.

Which is why, from one vantage, there are 7 distinct (but
isomorphic) Box-Kites in Sedenion space:  because we've 7 choices of
which Octonion to suppress!  6 vertices times 7 gives us the
\textit{42 Assessors} of our first ZD paper [3], a term we'll use
interchangeably with \textit{dyad} throughout.  We can, in fact, tug
on the network of interconnected lines ``wok-cooking'' style,
stirring things into and out of the hot oil in the center of the
Box-Kite.  (${\bf S}$ as "Stir-fry constant"?)  To find the
``Octonion copy'' labeling low indices on Box-Kite vertices where
the 5, say, is suppressed, trace the line containing it and the 4,
and ``rotate'': the 1 now goes from the left slope's midpoint to the
bottom right angle, to be replaced by the 4 while the 5 heads for
the middle, with CPO order (and hence, orientation of the line)
remaining unchanged.  Of the other 2 trips the 5 belongs to, only
one will preserve midpoint-to-angle orientation along the 6
o'clock-to-midnight vertical:  $(2, 5, 7)$, as one can check in an
instance.  (The two possibilities must orient oppositely when placed
along the same line, since one is Rule 1, the other Rule 2.)

From this point, everything is forced. This is obviously a procedure
that is trivial to automate, for any ``Octonion copy,'' regardless
of the ambient dimensionality the Box-Kite it underwrites might
float in. This simple insight will be the basis, in fact, of our
proof method, both in this paper and its sequel.  Another simple
insight will tell us how to find the \textit{high}-index term for
any vertex's dyad.  Two indices per vertex leaves 4 that are
suppressed: 0 (for the Reals), ${\bf G}$ and ${\bf S}$, and the XOR
(and also simple sum) of the latter two, which we'll shorthand ${\bf
X}$. These four clearly form a Quaternion copy -- one, in fact,
which has no involvement whatsoever in its containing Box-Kite's
zero-divisions.  Putting the index of the \textit{one} among these
which is \textit{itself} an L-unit center stage gives us the full
array of \textit{L-index sets} (trips composed of those indices of a
Sail's 3 vertices $< {\bf G}$) associated with the 4 Sails. Putting
in ${\bf G}$ or ${\bf X}$, then, must give us the full array of
\textit{U-index sets} (``U'' as in ``upper'').

Since each node belongs to 3 lines in PSL(2,7), the strut constant
belongs to 3 trips, each containing one term from the Rule 0 Zigzag
Sail's L-index set, and one from the Vent which resides opposite it
on the Box-Kite's octahedral frame.  In the Sedenions, three simple
rules govern interactions of the Vent and Zigzag dyads sharing a
strut. Writing the U- and L- index terms in upper and lower case
respectively, we can symbolize their dyads as (V, v) and (Z, z)
respectively.  The ``Three Viziers'' (derived as side-effects of
[VII], with one for each non-0 member of our ZD-free index set) read
as follows:

\begin{center}
VZ1:  $v \cdot z = V \cdot Z = {\bf S}$

\smallskip
VZ2:  $Z \cdot v = V \cdot z = {\bf G}$

\smallskip
VZ3:  $V \cdot v = z \cdot Z = {\bf X}$.
\end{center}

The First Vizier motivates the term \textit{strut constant}: for the
same pattern obtains for it, regardless of the strut being
investigated. The Second Vizier shows us that ${\bf G}$ connects
\textit{strut opposites}, always by Rule 1 logic.  But clearly, the
Third Vizier gives us the simplest way to answer any questions
concerning the relations between indices within a dyad: the L- and
U- indices of \textit{any} dyad belong to the same trip as ${\bf
X}$, with CPO ordering determined by whether or not the dyad belongs
to the Zigzag proper or the Vent opposite it.

Beyond the Sedenions, VZ2 is universally true, but VZ1 and VZ3 are
only so \textit{up to sign}:  e.g., the VZ1 L-trip for an arbitrary
strut can read $(z,c,S)$ in certain higher-dimensional contexts.
This is ultimately a side-effect of the same ``carrybit overflow''
that creates the phenomenon of most interest to us here, the
``missing box-kites'' in all $2^N$-ions, $N$ at least 5, for $S
> 8$ and not a power of 2. Correlated with such ZD-free structures are ``Type II''
box-kites with $S < 8$ (or, more generally, < {\bf G}/2),
indistinguishable from the standard ``Type I'' variety but for strut
orientations (with exactly 2 of a ``Type II'''s 3 struts always
being reversed:  see Appendix B).  Their ``twist products''
(operating similarly on parallel sides of each of the 3 orthogonal
squares or ``catamarans'' of a box-kite's orthogonal wire-frame, as
opposed to the 4 triangular ``sails'' which are our sole focus in
this monograph) let them act as middlemen between the normal and
ZD-free structures.  Our arguments here will make no use of such
``twist product'' subtleties (on which, see Theorem 6 in Part I and
the caveat that follows it, and the more developed remarks and
diagrams in [4]). Indeed, their phenomenology falls ``under the
radar'' of our Sail-based analysis:  strut-opposite Assessors, after
all, do \textit{not} mutually zero-divide.

Given our limited purposes here, therefore, our toolkit, once the
Viziers are dropped in it, is complete for all our later proofs. (We
must simply remember that invocations of VZ1 and VZ3 implicitly
concern \textit{sign-free} relations between Vent and Zigzag terms
-- that is, indices of XOR products only.)  What's left to do still:
get our hands messy with the plumbing, and then clean up with a last
grand construct. Let's start with the plumbing, and add some
notation. Label the Zigzag dyads with the letters A, B, C; label
their strut-opposite terms in the Vent F, E, D respectively. Specify
the diagonal lines containing all and only ZDs in any such dyad K as
$\verb|(K, /)|$ and $\verb|(K, \)|$ -- for $c \cdot (i_{K} + i_{k})$
and $c \cdot(i_{K} - i_{k})$ respectively, $c$ an arbitrary real
scalar. The twelve edges of the octahedral grid are so many pipes,
through which course the two-way streets of \textit{edge-currents}:
for the 3 edges of the Zigzag (and the 3 defining the opposite
Vent), currents joining arbitrary vertices M and N are called
\textit{negative}, since they have this form:

\begin{center}
$\verb|(M,/)| \cdot \verb|(N, \)|$ $= \verb|(M, \)|$ $\cdot$
$\verb|(N,/)|$ $ = 0$
\end{center}

Tracing the perimeter of the Zigzag with one's finger, performing ZD
products in natural sequence -- $\verb|(A, /)| \cdot $\verb|(B, \)|,
followed by the latter times $\verb|(C, /)|$, then this times
$\verb|(A, \)|$ and so forth -- one should quickly see how the
Zigzag's name was suggested.  Suppressing all letters, one is left
with just this cyclically repeating sequence:  \verb|/\/\/\|.

Currents along all 6 edges joining Zigzag and Vent dyads, on the
contrary, connect similarly sloping diagonals, hence are called
\textit{positive}, yielding the shorthand sequence \verb|///\\\| for
Trefoil sail traversals:

\begin{center}
$\verb|(Z,/)| \cdot \verb|(V, /)| = $ $\verb|(Z, \)| \cdot
\verb|(V,\)| = 0$
\end{center}

Consider the chain of ZD multiplications one can make along the
Zigzag, between A and B, then B and C, then C and A, for ${\bf S} =
4$.  The first term of this 6-cycle of zero products, once fully
expanded, is writable thus:

\begin{center}
$\verb|(A, /)| \cdot \verb|(B, \)|$ $= (i_{1} + i_{13}) \cdot (i_{2}
- i_{14}) = (i_{3} - i_{15} + i_{15} - i_{3}) = $ \newline
$\verb|(C, /)| - \verb|(C, /)| = \verb|(C, \)| - \verb|(C, \)| = 0$
\end{center}

We can readily see here where the notion of \textit{emanation}
arises:  traversing the edge between any two vertices in a Sail
yields a balance-pan pairing of oppositely signed instances of the
terms at the Sail's third vertex ... the 0 being, then, an instance
of ``balanced bookkeeping'' (whence the term ``Assessor,'' our
synonym for ``dyad'').  This suggests the spontaneous emanation of
particle/anti-particle pairings from the quantum vacuum, rather than
true ``emptiness.''

Finally, a side-effect of such ``Sail dynamics'' is this astonishing
phenomenon:  each Sail is an interlacing of 4 associative triplets.
For the Zigzag, these are the L-index $(a, b, c)$, plus the 3
U-index trips obtained by replacing \textit{all but one} of these
lowercase letters with their uppercase partners:  ergo, $(a, B, C)$;
$(A, b, C)$; $(A, B, c)$.  Ultimately this tells us that ZDs are
extreme \textit{preservers} of order, since they maintain
associativity in rigorous lock-step patterns, for all $2^N$-ions, no
matter how close to $\infty$ their $N$ might become.  Put another
way, the century-long aversion reaction experienced by virtually all
mathematicians faced with zero-divisors was profoundly misguided.

\section{Emanation Tables:  Conventions for Construction}
Theorem 7 guaranteed the simple structure of ETs:  because any
Assessor's uppercase index $i_{U}$ is strictly determined by ${\bf
G}$ and ${\bf S}$, once we are given these two values, the table
need only track interactions among the lowercase indices $i_{L}$.
This will only lead to ambiguities in the very place these are
meaningful:  in the recursive articulation of a boxes-within-boxes
tabulation of meta-fractal or Sky behaviors. In such cases, the
overlaying will be as rich in significance as the multiplicity of
sheets of a Riemann surface in complex analysis.

An ET does for ZD interactivity what a Cayley Table does for
abstract groups:  it makes things visible we otherwise could not see
-- and in a similar way.  Each Assessor's L-index is entered (in a
manner we'll soon specify) as a row (R) or column (C) value, with
XOR products (P values) among them being placed in the ``spreadsheet
cell'' (r,c) uniquely fixed by R and C.  We've noted such values
only get entered if P is the L-index of a legitimate emanation: that
is, the Assessor it represents mutually zero-divides (forms
\textit{DMZ}s with, for ``divisors making zero'') \textit{both} the
Assessors represented by the R and C labels of its cell. (As already
suggested, the natural use of the letters R, C, P here inspired
calling the study of NKS-like ``simple rules'' for cooking fractals
from their bit-strings \textit{recipe theory}.)

Four conventions are used in building ETs:  first, their labeling
scheme obeys the same nested-parentheses ordering we've already used
in designating Assessors A through F, with D, E, F the strut
opposites of A, B, C in reverse of the order just written.  The
L-indices, then, are entered as labels running across the top and
down the left.  The label of the lowest L-index is placed flush left
(abutting the ceiling), with the corresponding label of its strut
opposite being entered flush right (atop the floor).  As there will
always be ${\bf G} - 2$ (hence, an even number of) indices to enter,
repeating this procedure after each pair has been copied to
horizontal and vertical labels will completely exhaust them all.

Second convention:  As the point of an ET is to display all
legitimate DMZs, any cell whose R and C do \textit{not} mutually
zero-divide is left blank -- even if, in fact, there \textit{is} a
well-defined XOR value.  Hence, if R and C reference the same
Assessor, the XOR of their L-indices will be $0$; if they reference
strut opposites, the XOR will be ${\bf S}$.  But in both cases, the
cell (hence, the P value) is left blank.  All ``normal'' ETs, then,
will have both long diagonals populated by blank cells, while all
other cells are filled.

Third convention:  the two ZD diagonals associated with any Assessor
are not distinguished in the ET, although various protocols are
possible that would make doing so easy.  The reasons are parsimony
and redundancy:  rather than create longer, or twice as many,
entries, we assume both entries for the same Box-Kite edge will
contain the positive-sloping diagonal when the lower L-index appears
as the row label, else the negative-sloping diagonal when the higher
L-index appears first instead.  Such niceties won't concern us much
here: the key thing is that, in fact, all 24 filled cells of a
Box-Kite's ET entries can be mapped one-to-one to its ZD diagonals.
Recall, per Theorem 3, that both ZD diagonals of an Assessor form
DMZs with the same Assessor, according to the same edge-sign logic.
This leads us to the \dots

Fourth convention:  Although they are superfluous for many purposes,
edge signs provide critical information for others, and so are
indicated in all ETs provided here.  Each of a Box-Kite's 12 edges
conducts two currents -- one per ZD diagonal -- and does so
according to one or the other orientational option.  ZD diagonals
are conventionally inscribed so that the horizontal axis of their
Assessor plane is the L-indexed unit, while the vertical is the
U-indexed unit.  But even if this convention were reversed, the
diagonal leading from lower left quadrant to upper right would still
correspond to the state of synchrony implied by $\pm k(i_{L} +
i_{U})$:  for some Assessor U, we write $\verb|(U,/)|$.  Conversely,
the orthogonal diagonal indicative of anti-synchrony is written
$\verb|(U,\)|$.  If DMZs formed by the Assessors bounding an edge
are both of same kind, then we call the edge blue or notate it [+];
if Assessors U and V only form DMZs from oppositely oriented ZD
diagonals -- $\verb|(U,/)| \cdot \verb|(V,\)| = 0 \Leftrightarrow
\verb|(U,\)| \cdot \verb|(V,/)| = 0$ -- then we call the edge red or
notate it [-].  However, for ET purposes, since the red edges are
the most informative (all-red-edged Zigzags providing the stable
basis of Box-Kite structure, while all-red-edged DEF Vents play a
key role in twist-product interpreting -- a deep topic touched upon
in Part I, which won't concern us further here), we leave them
unmarked.  The six blue edges bounding the hexagonal view of the
Box-Kite, however, are preceded by an extra mark (best interpreted
as a dash, rather than a minus sign). This has the pragmatic
advantage that when zoomed, a large ET will have its entries with an
extra mark become unreadable in many software systems (e.g., one
sees only asterisks) -- and so we want the unmarked entries to be
those likely to be of most interest.

Since, given ${\bf X}$ (or, alternatively, ${\bf G}$ or N, and ${\bf
S}$), we can reconstruct a Box-Kite from just its Zigzag's L-index
trip, gleaning this information from an ET is worth explaining.  If
a given row contains the indices of any such Zigzag L-trips, they
will appear as the row label itself, plus two unmarked cell entries,
with the column label of the one appearing as the content of the
other. (If either cell in such a complementary set be marked with a
dash, then we are dealing with a DEF Vent index.) Each Zigzag L-trip
will also appear 3 times in an ET, once in each row whose label is
one of its indices, its 2 non-label indices appearing in un-dashed
cell entries each time.

Here is a readily interpreted emanation table.  Having $6 = 2^{3} -
2$ rows and columns, ${\bf G} = 8$, so ${\bf N} = 4$, making this a
Sedenion ET (encoding, thereby, a single Box-Kite).  And, since $2
\veebar 3 = 4 \veebar 5 = 6 \veebar 7 = 1$, the Strut Constant ${\bf
S} = 1$ as well.  A scan of the first row shows $6$ and $5$
unmarked, under headings $4$ and $7$ respectively; however, these
two labels appear as cell values which \textit{are} marked, making
these edges that connect Assessors in the D, E, F Vent.  In the
fourth row of entries, though, column labels $5$ and $3$ contain
cell values $3$ and $5$ respectively, both unmarked.  With their row
label $6$, then, these form the Zigzag L-index set $(3,6,5)$, which
hence must map to Assessors $(A, B, C)$.  Using the mirror-opposite
logic of the labeling scheme to determine strut opposites, it is
clear that the six row and column headings $(2, 4, 6, 7, 5, 3)$
correspond, in that order, to the Assessors $(F, D, B, E, C, A)$.
(The unmarked contents 6 and 5 in the first row, having labels
$(2,4)$ and $(2,7)$, thereby map to edges FD and FE, connecting DEF
Vent Assessors as claimed.) Finally, the long diagonals are all
empty: those cells in the diagonal beginning at the upper left all
have identical row and column labels; those in the mirror-opposite
slots, meanwhile, have labels which are strut-opposites.  By our
second convention, all these cells are left blank.
\begin{center}
\fbox{
\begin{tabular}{c | c |c |c |c |c | c}
  $\;$  & ${\bf 2}$   & ${\bf 4}$  & ${\bf 6}$  & ${\bf 7}$  & ${\bf 5}$  & ${\bf 3}$  \\ \hline
  \textbf{2}    & $\;$  & $6$  & $-4$ & $5$  & $-7$ & $\;$ \\ \hline
  \textbf{4}    & $6$   & $\;$ & $-2$ & $3$  & $\;$ & $-7$ \\ \hline
  \textbf{6}    & $-4$  & $-2$ & $\;$ & $\;$ & $3$  & $5$  \\ \hline
  \textbf{7}    & $5$   & $3$  & $\;$ & $\;$ & $-2$ & $-4$ \\ \hline
  \textbf{5}    & $-7$  & $\;$ & $3$  & $-2$ & $\;$ & $6$  \\ \hline
  \textbf{3}    & $\;$  & $-7$ & $5$  & $-4$ & $6$  & $\;$ \\
  \end{tabular}
  }
  \end{center}

Before beginning an in-depth study of emanation tables by type,
there is one general result that applies to them all -- and whose
proof will give us the chance to put the Three Viziers to good use.
While seemingly quite concrete, we will use it in roundabout ways to
simplify some otherwise quite complicated arguments, beginning with
next section's Theorem 9.  This Roundabout Theorem is our

\medskip

\noindent {\small Theorem 8.}  The number of filled cells in any
emanation table is a multiple of 24.

\medskip

\noindent \textit{Proof.}  Since 24 is the number of filled cells in
a Sedenion Box-Kite, this is equivalent to claiming that CDP
zero-divisors come in clusters no smaller than Box-Kites.  We have
already seen, in Theorem 5, that the existence of a DMZ implies the
3-Assessor system of a Sail, which further (as Theorem 7 spelled
out) entails a system of 4 interlocking trips:  the Sail's L-trip,
plus 3 trips comprising each L-trip index plus the U-indices of its
Assessor's 2 ``sailing partners.''  Since we have an ET, we have a
fixed ${\bf S}$ and fixed ${\bf G}$.  Hence, if we suppose our DMZ
corresponds to a Zigzag edge-current, we immediately can derive its
L-trip by Theorem 5, and all 3 Zigzag strut-opposites' L-indices by
VZ 1, and all 6 U-indices by VZ 3.  We then can test whether the
Trefoil Sails' edge-currents are all DMZs as follows.  As we wrote
in Theorem 7, $(u, v, w)$ maps to the Zigzag L-trip in CPO, but not
necessarily in $(a,b,c)$, order:  hence, $(u_{opp}, w_{opp}, v)$ is
an L-trip, and can be mapped to any of the Trefoils.  In other
words, given the Zigzag's 3-fold rotational symmetry, proving the
truth of the following arithmetical result proves the DMZ status of
\textit{all} Trefoil edges.  Yet we can avail ourselves of all 3
Zigzag U-trips in proving it.

\smallskip

\begin{center}
 $(w_{opp} - W_{opp})$

\underline{ $(u_{opp} + U_{opp})$}

$- V \;\; - v$

\underline{$ + v \;\; + V$}

$0$

\end{center}

\smallskip

The left bottom result is a given of the trip we started with.  The
result to its right is a three-step deduction from one of the Zigzag
U-trips:  use $(u_{opp}, w, v_{opp})$; Rule 2 gives $(u_{opp},
v_{opp} + G, w + G)$; the Second Vizier tells us this is $(u_{opp},
V, W_{opp})$; but the negative inner sign on the upper dyad reverses
the sign this trip implies, yielding $+ V$ for the answer.

The top results are derived similarly:  find which of the 4 Zigzag
trips underwrites the Vizier-derived ``harmonic'' which contains the
pair of terms being multiplied, and flip signs as necessary.  Hence,
the top left uses $(u, w_{opp}, v_{opp})$, then applies Rule 2 and
the Second Vizier to get $(- V)$, while the top right uses the
Zigzag L-trip itself:  $(u,v,w) \rightarrow (w + G, v, u + G)
\rightarrow (W_{opp}, v, U_{opp})$ -- which, multiplied by $(-1)$,
yields $(- v)$.{  } $\blacksquare$

\medskip

\noindent \textit{Remark.}  The implication that, regardless of how
large N grows, ZDs only increase in their interconnectedness, rather
than see their basic structures atrophy, flies in the face of a
century's intuition based on the Hurwitz Proof.  That there are no
standalone edge-currents, nor even standalone Sails, bespeaks an
astonishing (and hitherto quite unsuspected) stability in the realm
of ZDs.

\medskip

\noindent \textit{Corollary.}  An easy calculation makes it clear
that the maximum number of filled cells in any ET for any
$2^{N}$-ions is just the square of a row or column's length in
cells, minus twice the same number (to remove all the blanks in long
diagonals):  that is, $(2^{N-1} - 2)(2^{N-1} - 2) - 2 \cdot 2^{N-1}
+ 4 = (2^{2N-2} - 6 \cdot 2^{N-1} + 8) = (2^{N-1} - 4)(2^{N-1} - 2)
= 4 \cdot (2^{N-2} - 1)(2^{N-2} - 2)$.  By Roundabout, we now know
this number is divisible by 24, hence indicates an integer number of
Box-Kites. But two dozen into this number is just $(2^{N-2} -
1)(2^{N-2} - 2)/6$ -- the \textit{trip count} for the
$2^{N-2}$-ions! (See Section 2 of Part I.)  We have, then, the very
important Trip-Count Two-Step: \textit{The maximum number of
Box-Kites that can fill a $2^{N}$-ion ET = $Trip_{N-2}$.} We will
see just how important this corollary is next section.

\section{ETs for $N > 4$ and ${\bf S} \leq 7$}
One of the immediate corollaries of our CDP Rules for creating new
triplets from old ones is something we might call the Zero-Padding
Lemma:  if two $k$-bit-long bitstring representations of two
integers $R$ and $C$ being XORed are stuffed with the same number
$n$ of $0$s between bits $j$ and $j+1$, $0 \leq j \leq k$, their XOR
will, but for the extra $n$ bits of $0$s in the same positions, be
unchanged -- and so will the sign of the product $P$ of CDP-derived
imaginary units with these three bit-strings representing their
respective indices.

\textit{Examples.}  $(1,2,3) \rightarrow (2,4,6)
\rightarrow (4,8,12)$ [Add 1, then 2, 0s to the right of each
bitstring]

$(1,2,3) \rightarrow (1,4,5) \rightarrow (1,8,9)$ [Add 1, then 2, 0s
just before the rightmost bit in each bitstring]

$(3,4,7) \rightarrow (3,8,11) \rightarrow (3,16,19)$ [Add 1, then 2,
0s just after the leftmost bit in each bitstring]

\textit{Proof.}  Rule 1 will create a new unit of index ${\bf G} +
L$ from any unit of index $L < {\bf G}$, regardless of what power of
2 ${\bf G}$ might be.  Rule 2, meanwhile, uses any power of $2$
which exceeds all indices of the trip it would operate on, then adds
this ${\bf G}$ to two of the members of the trip, creating a new
trip with reversed orientation -- one of an infinite series of such,
differing only in the power of $2$ (hence, position of the leftmost
bit) used to construct them.  The lemma, then, is an obvious
restatement of the fundamental implications of the CDP Rules.

But creation of U-indices associated with L-indices in Assessor
dyads is the direct result of creating new triplets with G+S as
their middle term.  Hence, if we call the current generator $g$ and
that of the next higher $2^{N}$-ions ${\bf G}\;( = 2 \cdot g)$, then
if Assessors with L-indices $u$ and $v$ form DMZs in the Sedenions
for a given strut constant ${\bf S}$, their U-indices will increment
by $g$ in the Pathions, and zero division will remain unaffected. By
induction, the emanation table contents of the Sedenion (R,C,P)
entries will remain unchanged for all N, for all fixed ${\bf S} \leq
7$.  This leads us to

\medskip
\noindent {\small Theorem 9.}  All non-long-diagonal cell entries in
all ETs for all N, for all fixed ${\bf S} \leq 7$, will be filled.

\smallskip
\noindent \textit{Proof.} Keeping the same notation, the
$2^{N}$-ions will have $g$ more Assessors than their predecessors,
with indices ranging from $g$ itself to $2g - 1 \; ( = G - 1)$.
Consider first some arbitrary Zigzag Assessor with L-index $z < g$,
whose U-index is ${\bf G} + z \cdot {\bf S}$. (If it were a Vent
Assessor, or a Zigzag on a reversed-orientation strut in a ``Type
II'' box-kite, the second part of the expression would be reversed:
${\bf S} \cdot z$, per the First Vizier. This effects triplet
orientation, but not absolute value of the index, however, and it is
only the latter which matters at the moment.) Now consider the
Assessor whose L-index is the lowest of those new to the
$2^{N}$-ions, $g$. We know it is a Vent Assessor, in all Box-Kites
with  ${\bf S < g}$, of which there are 7 per each such ${\bf S}$ in
the Pathions, 35 in the 64-D $2^{6}$-ions, and so on: for it belongs
to the trip $({\bf S}, g, g + {\bf S})$ (Rule 1), so that its
U-index appears on its immediate left in the triplet $({\bf G} + g +
{\bf S}, g, {\bf G + S})$ (Rule 2 and last parentheses). Its
U-index, then, is ${\bf G} + (g \veebar {\bf S})$, or (recall Rule
1) just ${\bf G} + g + {\bf S}$. We claim these Assessors form DMZs;
or, writing out the arithmetic, that the following term-by-term
multiplication is true:

\smallskip
\begin{center}
$+ g + (G + g + S)$

\underline{$\;\;+ z + (G + z_ \cdot  S)\;\;$}
\smallskip

$- (G + g + z \cdot S) {  } - (z + g)$

\underline{$ + (z + g) {    } + (G + g + z \cdot S)$}

\smallskip
$0$\end{center}

\smallskip

Because one Assessor is assumed a Zigzag, while the other is proven
a Vent, the inner signs will be the same.  (Simple sign reversals,
akin to those involving our frequently invoked binary variable $sg$,
will let us generalize our proof to include the Vent-times-Vent case
later.)  Let's examine the terms one at a time, starting with the
bottom line.  Its left term is an obvious application of Rule 1, as
$z < g$, the latter being the Generator of the prior CDP level which
also contained $z$ as an L-index.  The term on bottom right we
derive as follows:  we know that $z$ and its U-index partner in the
$2^{N-1}$-ions belong to the triplet mediated by $g + {\bf S}$: $(z,
g + z \cdot {\bf S}, g + {\bf S})$. Supplementing this CPO
expression by adding ${\bf G}$ to the right-hand terms (Rule 2), we
get the triplet containing both multiplicands of the bottom-right
quantity:  $(z,  {\bf G} + g + {\bf S}, {\bf G} + g + z \cdot {\bf
S})$.  The multiplicands appear in this trip in their order of
application in forming the product; therefore, their resultant is a
plus-signed copy of the trip's third term, as shown above.

Moving to the left-hand term of the top line, what trip do the
multiplicands belong to?  Within the prior generation, Rule 1 tells
us that $z$'s strut opposite, $z \cdot {\bf S}$, multiplies $g$ on
the left to yield $g + z \cdot {\bf S}$.  Application of Rule 2 to
the terms $\neq g$ reverses order and gives us this: $({\bf G} + g +
z \cdot {\bf S}, g, {\bf G} + z \cdot {\bf S})$.  But what we've
written above is the product of multiplying the third and second
terms of the trip together, in CPO-reversed order; hence, the
negative sign is correct.  Finally, we get the negative of $(z + g)$
by similar tactics:  the term is the U-index of $z$'s strut-opposite
Assessor in the prior CDP generation, hence belongs to the trip with
this CPO expression:  $(g + {\bf S}, z + g, z \cdot {\bf S})$. Rule
2 gives us $({\bf G} + g + {\bf S}, {\bf G} + z \cdot {\bf S}, z +
g)$.  Hence, the product written above is properly signed.

Now, what effect does our initial assumption that $z$ is the L-index
of a Zigzag Assessor have on the argument?  The lower-left term is
obviously unaffected.  But the upper-left term, perhaps less
obviously, also is unchanged:  while it seems to depend on $z \cdot
{\bf S}$, in fact this is only used to define the L-index of $z$'s
strut opposite, which multiplies $g$ on the left to precisely the
same effect as $z$ itself, both being less than it.  The two terms
on the right, just as clearly, \textit{do} have their signs changed,
for in both, the order relations of L- and U- indices \textit{vis
\`a vis} ${\bf G + S}$ or ${\bf X}$ are necessarily invoked.  But
both signs on the right can be re-reversed to obtain the desired
result if we change the inner sign of the topmost expression --
which is to say, we have an effect analogous to that achieved in
earlier arguments by use of the binary variable $sg$, as claimed.

Since one CDP level's ${\bf G}$ is the $g$ of the next level up, the
above demonstration clearly obtains, by the obvious induction, for
all $2^{N}$-ions including and beyond the Pathions.  But what if one
or both L-indices in a candidate DMZ pairing \textit{exceed} $g$?

Rather than answer directly, we use the Roundabout Theorem of last
section.  Given a DMZ involving Assessors with L-indices $u < g$ and
$g$, we are assured a full Box-Kite exists with a Trefoil L-trip
$(u, g, g + u)$.  The remaining Assessors, being their strut
opposites, then have L-indices $u_{opp}, g + {\bf S}$, and $g + u
\cdot {\bf S}$.  As $u$ varies from $1$ to $7$, skipping ${\bf S} <
8$, zero-padding assures us that all DMZs from prior CDP generations
exist for higher N, for all L-indices $u, v < 8$.  Only those
Box-Kites created by zero-padding from prior-generation Box-Kites
(of which there can be but $1$ inherited per fixed ${\bf S}$ among
the $7$ found in the Pathions, for instance) will have all L-indices
$< g$.  For all others, the model shown with those having $g$ as an
L-index must obtain. Hence, only one strut will have L-indices $<
g$, the rest being comprised of some $w$ with L-index $\geq$ 8, the
others deriving their L-indices from the XOR of $w$ with the strut
just mentioned, or with ${\bf S}$.

But what will guarantee that any edge-currents will exist between
arbitrary Assessors with L-indices $u < g$ and $g + k, 0<k<g$, since
there is not even \textit{one} DMZ to be found among Assessors with
L-indices $\leq g$ in the candidate Box-Kite they would share? We
can now narrow the focus of our original question considerably, by
making use of the curious computational fact we called the
Trip-Count Two-Step.

In Part I's preliminary arguments concerning CDP, we showed that the
number of associative triplets in a given generation of
$2^{N}$-ions, or \textit{$Trip_{N}$}, can be derived from a simple
combinatoric formula.  Call the count of complete Box-Kites in an ET
$BK_{N,S}$. For ${\bf S} < 8$, $BK_{N,S} = Trip_{N-2}$, provided all
L-indices $g + k, 0 < k < g$, form DMZs in the candidate Box-Kites
implied. To begin an induction, let us consider a new construction
along familiar lines, which will provide us an easy way to
comprehend the Pathion trip-systems of all ${\bf S} < 8$. Beginning
with N = 5, we designate $Trip_{N-2}$ trips for each ${\bf S} < 8$
as type Rule 0, in the manner the singleton $2^{2}$-ion trip
$(1,2,3)$ was used in our introduction's "wok-cooking" discussion
(which Part I, Section 5, used as the basis of its ``slipcover
proofs''). But now, instead of putting the Octonions' ${\bf G = 4}$
in the center of the PSL(2,7) triangle, we put the Sedenions' $8$.

For consistency of examples, we continue to assume ${\bf S} = 1$, so
we'll begin with $(3,6,5)$, the Zigzag L-trip for ${\bf S = 1}$ in
the Sedenions, and also, by zero-padding, an L-trip Zigzag for 1 of
the 7 Box-Kites with ${\bf S = 1}$ among the Pathions. Extending
rays from the $(3,6,5)$ midpoints through the center creates Rule 1
trips which end in $11, 14, 13$: $(a,b,c)$ get sent to $(F, E, D)$
respectively.  The Rule 2 trips along the sides, in order of Zigzag
L-index inclusion, then correspond to Trefoil U-trips, all oriented
clockwise.  They read symbolically (literally) as follows: $EaD \;
(14, 3, 13); DbF \; (13,6,11); FcE
 \; (11,5,14)$.  We claim each of these 7 lines, when its nodes are
attached to their strut opposites, map 1-to-1 to an ${\bf S = 1}$
Pathion Box-Kite.  We have this as a given for the Rule 0 trip; we
need to \textit{explain} this for the Rule 1 trips (which Roundabout
already tells us are Box-Kites); and, we need to \textit{prove} it
for the Rule 2 trips that make the sides. (And, once we \textit{do}
prove it, and frame the suitable induction for all higher N, the
task which originally motivated us will be done: for these U-trips
house the Assessors with L-indices $> g$, whose candidate Box-Kites
don't include $g$.)

The Rule 1 trips, in all instances within this example, correspond
to Assessor L-indices $(a, d, e)$.  With $g = 8$ at $d$, the Third
Vizier tells us $c = 8 + {\bf S} = Sedenion \; {\bf X}$. $(a, b, c)$
thereby reads, within the Sedenions, as $(a, A, X)$. But in the
Pathions, all 3 terms are less than ${\bf G}$, hence can comprise an
L-index trip for a Sail -- and specifically, a Zigzag (else the
order of $A$ and ${\bf X}$ would be reversed). Similarly, the old
Sedenion $(f, F)$ are the new Pathion $(f, e)$, with the new trip
$(f,c,e)$ being the Third Vizier's way of saying $(f,X,F)$ from the
Sedenions' vantage.

For the Rule 2 trips, we prove one relation in one of them a DMZ,
which Roundabout tells us implies the whole Box-Kite, while symmetry
allows us to assume the same of the other two.  Consider, then, the
$aDE$ Trefoil U-trip, instantiated by $(3, 13, 14)$ in our example;
specifically, compute the product of the Assessors containing $a$
and $D = c + g$ as L-indices.  Their U-indices within the Pathions
must be $({\bf G} + a \veebar {\bf S}) = ({\bf G} + f)$, and $({\bf
G} + g + c \veebar {\bf S}) = ({\bf G} + g + d)$ respectively.  We
write their dyads when multiplying with opposite inner signs, as we
assume their DMZ is an edge in a Zigzag.  We claim the truth of this
arithmetic:

\smallskip
\begin{center}
$+ (c + g) - (G + g + d)$

\underline{$\;\;+ a   \;\;\;\;\;\;\;\;+\;\;\;\;\;  (G + f)\;$}
\smallskip

$+ (G + g + e) {  } - (b + g)$

\underline{$ + (b + g) {    } - (G + g + e)$}

\smallskip
$0$\end{center}

\smallskip

Bottom left:  $(a, b, c) \rightarrow (a, c + g, b + g)$ (Rule 2,
with $N=4$.)

Bottom right:  $(a,d,e) \rightarrow (a, g+e, g+d) \rightarrow (a,
{\bf G} + g + d, {\bf G} + g + e)$ (Rule 2 twice, $N=4$, then
$N=5$.) Upper dyad's inner sign reverses that of product.

Top left:  $(f,c,e) \rightarrow (e + g, c + g, f) \rightarrow ({\bf
G} + f, c + g, {\bf G} + e + g)$ (Rule 2 twice, $N=4$, then $N=5$.)

Top right:  $(f,d,b) \rightarrow (b + g, d + g, f) \rightarrow (b +
g, {\bf G} + f, {\bf G} + g + d)$  (Rule 2 twice, $N=4$, then
$N-5$.) Upper dyad's inner sign reverses that of product.

A similar brief exercise with either DMZ formed with the emanated
Assessor will show it, too, has a negative inner sign with respect
to a positive in its DMZ partner. Two negative edge-signs in one
Sail means Zigzag (means \textit{three} negative edge-signs, in
fact). Our proof up through the Pathions is complete; we need only
indicate the existence of a constructive mechanism for pursuing this
same strategy as $N$ grows arbitrarily large.

Consider now the same PSL(2,7) triangle, but in its center put a 16
($= {\bf g} = {\bf G/2}$ for the 64-D Chingons, after the 64
Hexagrams of the \textit{I Ching}, to give them a name).  Then, put
all 7 of the Pathions' ${\bf S = 1}$ Zigzag L-trips into the Rule 0
circle.  One gets $3 \cdot 7 = 21$ Rule 2 Zigzag L-trips, and the 10
integers $< g$ found in them and the 7 Rule 0 Zigzag L-trips implies
there are 10 Rule 1 Trefoil L-trips, each associated with a distinct
Box-Kite.  But that would make for $7 + 21 + 10 = 38$ Zigzag
L-trips, when we know there can only be $35$. The extra 3 indicate
there's some double-duty occurring: specifically, 3 of the Rule 1
Trefoil L-trips in fact designate not the standard $(a,d,e)$, but
$(f,d,b)$, with $d = g = 16$ in each instance.  When $(5,14,11)$ is
fed into our ``trip machine'' as Rule 0 circle, both $(11,16,27)$
and $(14,16,30)$ map to $(f,d,b)$ trips tied to Rule 0 Zigzag
L-indices $(10,27,17)$ and $(15,30,17)$, whose $(a,d,e)$ trips
appear as rays on triangles for $(3,10,9)$ and $(3,13,14)$
respectively.  $(11,16,27)$ also shows as an $(f,d,b)$ with Rule 0
trip $(6,11,13)$. (Readers are encouraged to use the code in the
appendix to [4], to generate ETs for low ${\bf S}$ and $N$.
Trip-machining details for our ${\bf S = 1}$ example are in Appendix
A.) For $N = 7$, use the 35 just-derived ${\bf S} = 1$ L-trips as
Rule 0 circles with a central $32$, and so on.{ } $\blacksquare$

\section{The Number Hub Theorem (${\bf S} = 2^{N - 2}$) for $2^{N}$-ions}

Given the lengths required to prove the fullness of ETs for ${\bf S
< 8}$, it might be surprising to realize that the infinite number of
cases for ${\bf S} = 2^{N - 2}$ for all $2^{N}$-ions are so simple
to handle that they almost prove themselves.  Yet the proof of this
Number Hub Theorem, while technically trivial, has far-reaching
implications.

\medskip
\noindent {\small Theorem 10.}  For all $2^{N}$-ions with ZDs ($N >
3$), and ${\bf S} = g = {\bf G/2}$, all non-long-diagonal entries in
the emanation table are filled; more, each such filled cell in the
ET's upper left quadrant is unmarked (indeed, indicates an
edge-current in a Zigzag); further, the row, column, and cell
entries are isomorphic to those found in an unsigned, CDP-generated,
multiplication table for the $2^{N-2}$-ions; finally, the
$Trip_{N-2}$ Zigzag L-index sets which underwrite its Box-Kites are
precisely all and only those trips contained in said $2^{N-2}$-ions,
the ET effectively serving as their high-level atlas.

\medskip

\noindent \textit{Proof.} As the largest L-index of any Assessor is
$2g - 1$, and each ${\bf S}$ in the ETs in question is precisely
$g$, then the row (column) labels will ascend from $1$ to $g - 1$ in
simple increments from top to bottom (left to right) in the upper
left quadrant, making its square of filled cells isomorphic to
unsigned entries in the corresponding $2^{N-2}$-ion multiplication
table. Also, all these filled cells of the ET will only contain XORs
of indices $< g$. Hence, all and only L-index trips will have the
edges of their (necessarily Zigzag) Sails residing in said quadrant.
All non-long-diagonal cells in the ET are meanwhile filled, since
all candidate Assessors have form $M = (m, {\bf G} + g + m)$, and
for any CPO triplet $(a,b,c)$ whose row and column labels plus cell
entry are contained in the upper left quadrant, it is easy to show
that the following arithmetic is true:

\begin{center}
$+ b - (G + g + b)$

\underline{$\;\;+ a   \; + \;  (G + g + a)\;\;$}
\smallskip

$+ (G + g + c) {  } - c$

\underline{$ + c {    } - (G + g + c)$}

\smallskip
$0$\end{center}

\smallskip

Therefore, the $Trip_{N-2}$ Box-Kites, the Zigzag L-index set of
each of which is one of the $Trip_{N-2}$ trips contained in the
$2^{N-2}$-ions, all have this simple form:

\begin{center}
$(a,b,c,d,e,f) = (a,b,c, g+c, g+b, g+a) \;\;\; \blacksquare$
\end{center}

\medskip

\noindent \textit{Remarks.}  As will become ever more evident,
powers of $2$ -- which is to say, singleton $1$-bits in indefinitely
long binary bitstrings -- play a role in ZD number theory most
readily analogized to that of primes in traditional studies. And
while integer triples (from Pythagoras to Fermat) play a central
role in prime-factor-based traditional studies, all XOR triplets at
two CDP generations' remove from the power of $2$ in question are
collected by its ET in this new approach. All other integers
sufficiently large (meaning $> 8$) are meanwhile associated with
fractal signatures, to each of which is linked a unique
infinite-dimensional space spanned by ZD diagonals.  But can such a
vantage truly be called Number Theory at all?  We say indeed it can:
that it is, in fact, the ``new kind of number theory'' that must
accompany Stephen Wolfram's New Kind of Science. In his massive 2002
book, he tells us that, common wisdom to the contrary, complex
behavior can be derived from the simplest arithmetical behavior. The
obstacle to seeing this resides in the common wisdom itself [5, p.
116]:

\begin{quote}
$\cdots$ traditional mathematics makes a fundamental idealization:
it assumes that numbers are elementary objects whose only relevant
attribute is their size.  But in a computer, numbers are not
elementary objects.  Instead, they must be represented explicitly,
typically by giving a sequence of digits.
\end{quote}

But that ultimately implies strings of $0$'s and $1$'s, where the
matter of importance becomes which places in the string are held,
and which are vacant:  the original meaning of our decimal
notation's sense of itself as placeholder arithmetic.  The study of
zero \textit{divisors} -- placeholder substructures -- then becomes
the natural way to investigate the composite characteristics of
Numbers \textit{qua} bitstrings.  When we discover, in what follows,
that \textit{composite} integers (meaning those requiring multiple
bits to be represented) are inherently linked, when seen as
strut-constant bit-strings, with infinite-dimensional meta-fractals,
the continuation of the quote on the following page should ring
true:

\begin{quote}
In traditional mathematics, the details of how operations performed
on numbers affect sequences of digits are usually considered quite
irrelevant.  But $\cdots$ precisely by looking at such details, we
will be able to see more clearly how complexity develops in systems
based on numbers.
\end{quote}

\section{ The Sand Mandala Flip-Book ($8 < {\bf S} < 16, \; N = 5)$}

In the first concrete exploration of ZD phenomenology beyond the
Sedenions [6, pp. 13-19], a startling set of patterns were
discovered in the ETs for values of ${\bf S}$ beyond the ``Bott
limit'':  that is, for $8 < {\bf S} < 16$ (the upper bound being the
${\bf G}$ of the 32-D Pathions), the filled cells sufficed to define
not $7$, but only $3$, Box-Kites for $N = 5$; more, the primary
geometric figures in each such ET transformed into each other with
each integer increment of ${\bf S}$, in a manner exactly reminiscent
of the flip-books which anticipated cartoon animation.  While these
seemed perplexing in mid-2002 when they were found, their logic is
in fact profoundly simple.

First, each such ET's ${\bf S}$ is just the ${\bf X}$ of one already
seen in the Sedenions.  We continue our convention of using $g$ to
indicate the ${\bf G}$ of the prior CDP generation, employ $s$ for
said generation's ${\bf S}$, and reference all prior Assessor
indices by suffixing their letters with asterisks.  Then, since
${\bf S} = g + s$, the trip $(s, g, g + s)$ mandates, by the First
Vizier (whose signed version we invoke due to the direct derivation
from the Sedenions), that $g$ must belong to the Zigzag Sail if it's
to be an Assessor L-index at all.

Note that this is not a truly legitimate argument, as we'll see
shortly, albeit the results are correct, as shown by other means in
[6]: this is because ``Type II'' box-kites first emerge in this
current context -- but are \textit{not} among the 3 x 7
``flip-book'' denizens of immediate interest.  We will assume, for
simplicity of presentation, that the First Vizier \textit{does}
obtain here:  \textit{proving} that it does, however, requires a
background argument concerning ``Type II'' box-kites:  their ${\bf
S}$ values must be \textit{less} than $g$, hence none of our
flip-book candidates can qualify.  (But they are just as numerous as
the flip-book box-kites, there being 3 for each of the seven values
of ${\bf S} < g$.  For their listing, and theoretical framing of
``Type II'' phenomenology, see Appendix B.)

We will content ourselves here with giving this as an empirical
result, and assume, therefore, the validity of the signed version of
the First Vizier in the case at hand.  Based on this assumption, we
can further claim that the Sedenion Vent L-indices, $f*, e*, d*$,
must also be associated with Zigzag Assessors.  By an argument
exactly akin to that of last section, we then have 3 candidate
Box-Kites to consider:  since the 3 Vent L-indices are all less than
$g$, they must be mapped to the 3 Assessors A, $g = 8$ must adhere
to B (and $s = 1$ to E), while the L-indices of the C Assessors
associated with $f*, e*, d*$ must be A*, B*, C* respectively.  The
proof is easy:  taking the new A, C Assessors $ = (f*, {\bf G} + g +
a*)$ and $(g, {\bf G} + s)$ in that order as readily generalizable
representatives, we do the arithmetic.

\smallskip

\begin{center}
$+ g - ({\bf G} + s)$

\underline{$\;\;+ f*   \; + \;  ({\bf G} + g + a*)\;\;$}
\smallskip

$+ ({\bf G} + a*) \;\;\;\; - (f* + g)$

\underline{$ + (f* + g)  - ({\bf G} + a*)$}

\smallskip
$0$\end{center}

\smallskip

The bottom left is just Rule 1.  For the bottom right, start with
the First Vizier:  $(f*, a*, s) \rightarrow (f*, {\bf G} + s, {\bf
G} + a*) \rightarrow (f*) \cdot ( - ({\bf G} + s) ) = - ({\bf G} +
a*)$. The top left is derived thus: $(a*, \; g, \; g + a*)
\rightarrow (g, \; {\bf G} + a*, \; {\bf G} + g + a*) \rightarrow
({\bf G} + g + a*) \cdot g = + ({\bf G} + a*).$  Finally, $(a*, s,
f*) \rightarrow (g + a*, g + f*, s) \rightarrow ({\bf G} + g + a*,
{\bf G} + s, g + f*)$, but the negative inner sign of the top dyad
reverses sign as shown.

The 3 Box-Kites thus derived are the only among the 7 candidates to
be viable:  for the Zigzag L-index of the ${\bf S = 1}$ Sedenion
Box-Kite does \textit{not} underwrite a Sail; hence, by what lawyers
would call a ``fruit of the poisoned tree'' argument, neither do the
3 U-trips associated with the same failed Zigzag. Using A* and B*,
then invoking the Roundabout Theorem, we see this readily:

\smallskip
\begin{center}
$+ b* + ({\bf G} + g + e*)$

\underline{$\;\;+ a*   \; + \;  ({\bf G} + g + f*)\;\;$}
\smallskip

$- ({\bf G} + g + d*) {  } - c*$

\underline{$ + c* {    } - ({\bf G} + g + d*)$}

\smallskip
NOT ZERO (only c*'s cancel)

\end{center}

\smallskip

With the appending of two successive bits to the left, the
bottom-left and top-right products are identical to those obtaining
without the $({\bf G} + g)$ being included.  Similarly, the top-left
product uses Rule 2 twice, to similar effect, but with $({\bf G} +
g)$ included in the outcome:  since $(f*,d*,b*)$ is CPO, we then get
$-({\bf G} + g + d*)$.  For the top-right result, meanwhile, the two
high bits induce a double reversal, then are killed by XOR, leaving
the product the same as if they hadn't been there: $(f*, c*, e*)
\rightarrow (g + e*, c*, g + f*) \rightarrow ({\bf G} + g + f*, c*,
{\bf G} + g + e*)$, hence $-c*$.  We have an argument reminiscent of
Theorem 2:  depending on the inner sign of the upper dyad, one pair
of products cancels or the other, but not both.

We see, then, that the construction given without explanation at the
end of Part I is correct.  The arguments given there concerning the
vital relationship of a Box-Kite's non-ZD structures to semiotic
modeling suggest that this ``offing'' (to use the appropriately
binary slang linked to Mafia hitmen) of a Zigzag's 4 triplets should
have a similarly significant role to play in such modeling.  This
has bearing not just on semiotic, but physical models, since the key
dynamic fact implicit in the Zigzag L- and U- trips (or just Z-trips
henceforth) is their similarity of orientation:  since $(a,b,c);
(a,B,C); (A,b,C); (A,B,c)$ are all CPO as written, we are
effectively allowed to do pairwise swaps of upper- and lower- case
lettering among them without inducing anything a physicist might
deem observable (e.g., a $180\,^{\circ}$ reversal or ``spin
quantum''). This condition of \textit{trip sync} breaks down as soon
as we attempt to allow similar swapping between Z-trips and their
Trefoil compatriots: in particular, those 2 which don't share an
Assessor with the Zigzag.  The toy model of [7] would use these
features to designate the basis of a ``Creation Pressure'' that
leads to the output of the string theorist's $E_{8} \times E_{8}$
symmetry. This symmetry, as discussed there, breaks in the standard
models when one of the primordial $E_{8}$'s decays into an $E_{6}$
-- which has $72$ roots to parallel the $72$ filled cells of our
Sand Mandalas.  For present purposes, the key aspect of this
correspondence is that, in ZD theory at least, the explosion of a
singleton Box-Kite into a Sand Mandalic trinity throws the
off-switch on the source of the dynamics:  the Z-trips which
underwrite trip sync no longer even underwrite Box-Kites.  The whole
scenario suggests nothing so much as those boxes which, when opened
by pushing an external lever, emit an arm which pulls up on the same
lever, forcing the box to close and the arm to return to its hiding
place inside it.

Let's turn now to the ET graphics of the flip-book sequence, so
suggestive of cellular automata.  For each of the 7 ET's in
question, all labels $< g$ are monotonically increasing, since ${\bf
S}$, and hence their strut opposites, exceed them all.  But the only
filled (but for long-diagonal crossings) rows and columns will be
those with labels equal to ${\bf S - g} = s$ and its strut-opposite
$g$, for these L-indices reside at E and B respectively in all 3
Box-Kites in the ensemble, hence either dyad containing one of them
makes DMZs within each of the trio's $(a,d,e)$ and $(f,d,b)$ Sails,
filling all 12 ($= 2^{4} - 2$, minus $2$ for diagonals) fillable
cells in each row or column tagged with these Assessors' label.
Thus, as $s$ is incremented, two parallel sets of perpendicular
lines of ET cells start off defining a square missing its corners,
then these parallels move in unit increments toward each other,
until they form a 2-ply crossbar once $s = 7$ (${\bf S = 15}$).  24
cells each have row label R or column label C $= s$; 24 reside in
lines with label $= g$; and 24 more have their contents P $= s$ or
$g$: these last have an orderliness that is less obvious, but by the
last ET in the flip-book, they have arrayed themselves to form the
edges of a diamond, orthogonal to the long diagonals and meeting up
with the crossbar at its four corners, with $s = 7$ values filling
the upward-pointing edges, and $g = 8$'s those sloping down.

The graphics for the flip-book first appeared in [6, p. 15]; they
were recycled on p. 13 of [8]; larger, easily-read versions of these
ETs were then included (along with numerous other Chingon-based
flip-books and other graphics we'll discuss later) as Slides 25-31
of the Powerpoint presentation comprising [1], delivered at Wolfram
Science's June 15-18, 2006, NKS conference in Washington, D.C. All
three of these resources are available online, and the reader is
especially encouraged to explore the last, whose 78 slides can be
thought of as the visual accompaniment to this monograph.
(Henceforth, references to numbered Slides will be to those
contained and indexed in it.)

\section{64-D Spectrography:  3 Ingredients for ``Recipe
Theory''}

In a manner clearly related to Bott periodicity, strut constants
fall into types demarcated by multiples of 8.  But unlike the
familiar modulo 8 categorization of types demonstrated, perhaps most
familiarly, in the Clifford algebras of various dimensions, the
situation with zero-divisors concerns not typology (which keeps
producing new patterns at all dimensions), but granularity.  As we
shall see, emanation tables for ${\bf S} > 8$ (and not a power of
$2$), aside from diagonally aligned cells in otherwise empty
stretches, display checkerboard layouts of parallel and
perpendicular \textit{near-solid lines} (NSLs), whose cells all have
emanations save for a pair of long-diagonal crossings, and whose
visual rhythms are strictly governed by ${\bf S}$ and $8$ or the
latter's higher multiples.

The rule we found in the 32-D Pathions for the Sand Mandalas
indicates that the basic pattern (and $BK_{5,\;{\bf S}}$ for $8 <
{\bf S} < 16$) is ``essentially the same'' for all of them.  We put
the qualifying phrase in quotes, as it is an open question at this
point what features, residing at what depth, \textit{are} indeed
``the same,'' and which are different.  For the moment, we will
invoke the term \textit{spectrographic equivalence} as a sort of
promissory note, hoping to stuff ever more elements into its
grab-bag of properties, beginning with two.  First is something at
once intuitively obvious but not readily proven.  (We will include a
corollary to a later theorem when we have done so).  Since the first
8 possible strut-constant values all display maximally-filled ETs,
and since anomalies displayed by higher values are strictly
side-effects of bits to the left of the 8-bit (which are, of course,
its multiples), it is natural to assume that any recursive induction
upon simpler forms will echo this ``octave'' structure: that each
time ${\bf S}$ passes a new multiple of $8$, it participates in a
new type. (As with the Sand Mandalas, we will see this means that
$BK_{N,\;{\bf S}}$ for the new 7- or 8-element spectral band of new
forms will differ from that found in its predecessor band.) This
will lead, in the most clear-cut cases -- $\;{\bf S} = 15$, or a
multiple of $8$ not a power of $2$, say -- to grids composed of 8 x
8 boxes some or all of whose borders are NSLs.

How we determine which cases are clear-cut, meanwhile, and why and
how we might want or need to privilege them, leads to our second
property to include up-front in our grab-bag.  In a manner
reminiscent of the various tricks -- like minors and cofactors --
used in classical matrix theory to prove two matrices are
equivalent, we can transform members of a spectral band into each
other by certain formal methods of hand-waving. With the Sand
Mandalas, for instance, we could replace concrete indices in the row
and column labels with abstract designations referencing the
$(a,b,c)$ values of each of their 3 Box-Kites, listed in one of a
number of predetermined orders:  by least-first CPO ordering of such
$(a,b,c)$ triplets, in a sequence determined by the Zigzag L-trip of
the Sedenion Box-Kite we can derive them from, for instance (which
is equivalent to the 3 sand-mandalic Box-Kites' $d$ values, as we've
seen).

Since which cells are filled is strictly determined by ${\bf S}$ and
${\bf G}$, such designations eliminate all individuality among the
ETs in question.  Hence, if certain display features of one of them
seem convenient, we can convert its ``tone row'' of indices
populating its row and column labels into an abstract layout,
governed by which index is associated with which Assessor, in the
manner sketched last paragraph.  We could then use \textit{this}
layout as the template for re-writes of all other ETs in the same
spectral band, knowing that results obtained using the specific
instantiation of the band could thereby be converted into exactly
analogous ones for the other band-members.

We will, in fact, implicitly adopt this tactic by using ${\bf S =
1}$ as an exemplary ``for instance'' in numerous arguments, while
employing the highest-valued ${\bf S}$ found among the Sand
Mandalas, $15$, to simplify the visualizing (and calculating) of
recursive pattern creation for fixed-${\bf S}$, growing $N$
sequences.  (${\bf S = 15}$ is chosen because it has all its low
bits filled, hence all XORs are derived by simple subtraction,
leaving carrybit overflow to show itself only in what matters most
to us:  the turning off of $4$ candidate Box-Kites in the Pathions,
and -- as we will show two sections hence -- $16$ in the Chingons,
and $4^{N-4}$ in all higher $2^{N}$-ions.) Where we termed, for
reasons already explained, the fixed-$N$, growing ${\bf S}$
sequences \textit{flip-books}, we designate these new displays (for
reasons we'll justify shortly) \textit{balloon-rides}.

While there is but one abstract type for the Sedenions, with one
Box-Kite for each of the 7 possible ${\bf S}$ values, a second
spectral band emerges in the Pathions to include the Sand Mandalas,
and two more are added for the 64-D Chingons.  By induction from the
universally shared first band for all $N > 3$, where there are
$Trip_{N-2}$ Box-Kites in each ET, for each ${\bf S \leq 8}$, the
first new spectrographic addition includes the upper multiple of $8$
that bounds it, since it is not a power of $2$: $16 < {\bf S} \leq
24$. The second new range, though, is bounded by ${\bf G}$, hence
does not include it, as it is tautologically a power of $2$ (which
powers, as we saw two sections ago, comply with a type all their
own, with the same Box-Kite-count formula as for the lowest spectral
band): $24 < {\bf S} < 32$.

Each of these two new bands displays a distinctive feature which
underwrites one of the three key ingredients for the \textit{recipe
theory} we are ultimately aiming for.  We call these, for ${\bf S}$
ascending, \textit{(s,g)-modularity} and \textit{hide/fill
involution} respectively.  The third key ingredient, meanwhile,
resides in the band that first emerges in the Pathions -- and whose
echo in the Chingons has recapitulative features sufficiently rich
as to merit the name of \textit{recursivity}. We will be devoting
Part III's first post-introductory section to a thorough treatment
of the simplest instance of this third ingredient, showing how to
ascend into the meta-fractal we call the Whorfian Sky (named for the
great theorist of linguistics, Benjamin Lee Whorf, whose last-ever
lecture on ``Language, mind and reality'' described the layering of
meaning in language in a manner strongly suggesting something akin
to it). Among many visionary passages in his descriptions of a
future cross-disciplinary science, the following seems most apt to
serve as the lead-in quote for the third and final sweep of our
argument [9]:

\begin{quote}
Patterns form wholes, akin to the Gestalten of psychology, which are
embraced in larger wholes in continual progression.  Thus the cosmic
picture has a serial or hierarchical character, that of a
progression of planes or levels.  Lacking recognition of such serial
order, different sciences chop segments, as it were, out of the
world, segments which perhaps cut across the direction of the
natural levels, or stop short when, upon reaching a major change of
level, the phenomena become of quite different type, or pass out of
the ken of the older observational methods.  But $\cdots$ the facts
of the linguistic domain compel recognition of serial planes, each
explicitly given by an order of patterning observed.  It is as if,
looking at a wall covered with fine tracery of lacelike design, we
found that this tracery served as the ground for a bolder pattern,
yet still delicate, of tiny flowers, and that upon becoming aware of
this floral expanse we saw that multitudes of gaps in it made
another pattern like scrollwork, and that groups of scrolls made
letters, the letters if followed in a proper sequence made words,
the words were aligned in columns which listed and classified
entities, and so on in continual cross-patterning until we found
this wall to be -- a great book of wisdom! [10, p. 248]
\end{quote}

\section*{Appendix A:  Genealogy of ${\bf S = 1}$ Box-Kites}

$N=4$:  Unique Quaternion L-index set $(1,2,3)$ fed as Rule~0 circle
into PSL(2,7) with central $g = 4$, yielding 7 Octonions trips, each
with a different ${\bf S}$.  For ${\bf S = 1}$, have $(3,6,5)$,
which becomes singleton Rule~0 for next level.

\medskip

\noindent $N=5$:  $(3,6,5)$ fed as Rule~0 circle into PSL(2,7) with
central $g = 8$ yields 3 Rule~2 L-trips as triangle's sides, which
(upon affixing their strut opposites as L-indices) generate (along
with zero-padded $(3,6,5)$ ) 4 Box-Kites with ${\bf X = G + 1 =
17}$. Triangle's medians become $(a,d,e)$ Trefoil L-index sets of 3
Rule~1 ${\bf S = 1}$ Box-Kites, making 7 in all.  These Zigzag
L-index sets become Rule~0 trips for the next level, and are:

\begin{center}
Rule 0:$\;\;(3,6,5)$

Rule 1:$\;\;(3,10,9); \; (6,15,9); \; (5,12,9)$

Rule 2:$\;\;(3,13,14); \; (6,11,13); \; (5,14,11)$

\end{center}

\medskip

\noindent $N=6$:  The 7 $N=5$ Zigzag L-index sets just listed are
fed as Rule~0 circles into PSL(2,7) triangles with central $g = 16$,
and are Zigzag L-index sets in their own right for Box-Kites with
${\bf X = G + 1 = 33}$.

10 Rule~1 medians, 3 redundant (as they generate $(f,d,b)$'s where
$(a,d,e)$'s are also given:  $(14,16,30$)* and $(11,16,27)$** in
$(5,14,11)$'s triangle, the latter also in $(6,11,13)$'s).  They are
associated with these 7 Zigzag L-index sets:

\begin{center}
$(3,18,17); \; (5,20,17); \; (6,23,17); \; (9,24,17);$

$(10,27,17)*; \; (12,29,17); \; (15,30,17)**$

\end{center}

Rule~2 sides:  3 per each Rule~0 trip, as follows:

\begin{center}
$(3,6,5) \rightarrow (3,21,22); \; (6,19,21); \; (5,22,19)$

$(3,10,9) \rightarrow (3,25,26); \; (10,19,25); \; (9,26,19)$

$(6,15,9) \rightarrow (6,25,31); \; (15,22,25); \; (9,31,22)$

$(5,12,9) \rightarrow (5,25,28); \; (12,21,25); \; (9,28,21)$

$(3,13,14) \rightarrow (3,30,29); \; (13,19,30); \; (14,29,19)$

$(6,11,13) \rightarrow (6,29,27); \; (11,22,29); \; (13,27,22)$

$(5,14,11) \rightarrow (5,27,30); \; (14,21,27); \; (11,30,21)$

\end{center}

\medskip

\noindent $N=7$:  Feed the just-listed 35 Zigzag L-index sets to
PSL(2,7)'s with $g = 32$, as Rule~0 circles, thereby generating the
155 ${\bf S = 1}$ Zigzags found in the $2^{7}$-ions, or Routions --
named for the site of the Internet Bubble's once-famed
``Massachusetts Miracle,'' Route 128 -- and so on.

\section*{Appendix B:  A Brief Intro to ``Type II'' Box-Kites}

The recursive generation of Zigzag L-sets just presented calls for
some close attention when the box-kites involved are Type II, since
they then have the diagonals of their PSL(2,7) triangles oriented
differently:  instead of all 3 leading from midpoints of the Rule 2
sides to the corners, only 1 of these will preserve orientation for
a Type II (with the other two having ``reversed VZ1'' rules in
evidence).  We first give a construction for producing all the Type
II box-kites in the Pathions, and then indicate the manner in which
their workings are intimately connected with the phenomenology of
twist products broached in Part I's Theorem 6.

The construction was presented with different framing in [8], where
we deployed a ``stereo Fano'' representation using side-by-side
triangles, the left being a proper PSL(2,7).  Within the Pathions,
there are 7 distinct box-kites for each ${\bf S}$ \textit{except}
for the ``flip-book'' trios, one for each ${\bf S}> 8$.  And for
${\bf S} = 8$ exactly, we saw in our discussion of the Number Hub
Theorem that we can build all 7 by placing 8 in the center of the
standard Fano (what we'll call PSL(2,7) henceforth), then taking the
Zigzag L-trip for each Sedenion ${\bf S}$ and placing its units at
the sides' midpoints, in the usual CPO order (in left, right, and
bottom sides respectively).  Each of these 7 lines then generates a
new box-kite in the \textit{Pathions} for the \textit{Sedenion}
${\bf S}$ in question.

If we re-inscribe the starter-kit L-trip, but change ${\bf G}$ to
the Pathion's 16, applying VZ2 gives us new U-index terms, but the
L-index terms for all 6 Assessors remain the same as for the
Sedenion box-kite:  we call this ``Rule 0'' instance the
\textit{zero-padded} box-kite (or just ZP) for the ${\bf S}$ value
in question.

If we take the 3 ``Rule 1'' triplets along the struts, and place
them not at the A, B, C positions of our new Pathion box-kites, but
instead at A, D, E (with 8 always winding up at D), we generate 3
more standard (Type I) box-kites.  For ${\bf S} = 1$, the Sedenion
Zigzag L-trip is just $(3, 6, 5)$, and each of its units becomes the
low-index `A' for a new Pathion box-kite, with L-indices written in
``nested parentheses'' order (that is, A, B, C, D, E, F) as follows:
$(3, 10, 9, 8, 11, 2)$; $(6, 15, 9, 8, 14, 7)$; $(5, 12, 9, 8, 13,
4)$.

But if we take the 3 ``Rule 2'' triplets along the edges, mapping
the Zigzag unit at the center of each to the low-index `A' of a new
box-kite, the 8 doesn't show at any Assessor, and two of the three
struts will have orientations reversed.  These ``Type II''
box-kites, again for ${\bf S} = 1$, written per the same convention
just used for their 8-bearing siblings, read like this:  $(3, 13,
14, 15, 12, 2)$; $(6, 11, 13, 12, 10, 7)$; $(5, 14, 11, 10, 15, 4)$.
Since the A and F low-index terms are the same as in the same-${\bf
S}$ \textit{Sedenion} box-kite, the strut they make obviously has
the standard orientation.  (But note that there is nothing essential
about the (A,F) strut here: the placing of the lowest-indexed unit
of the Zigzag L-trip at A is a convenient convention, and its
employment in the Pathions suffices to induce this effect; however,
it no longer suffices in higher dimensions, where ${\bf S}$ can
exceed 8 yet still be less than $g$.)

That their being Type II is an immediate side-effect of ``Rule 2''
in this method of deriving them should be obvious.  What is less
obvious is their special relationship with twist products.  Here, we
review some of the basics:  in the Sedenions, whenever two Assessors
bound an edge, we can swap a pair of corresponding terms (either L-
or U- indices) and then switch the sign joining the L- and U-
indices in the resultant pairing, and get an Assessor in another
box-kite as a result.  Such ``twist products,'' then, reverse the
edge-sign of a given line of ZDs as we move between containing
box-kites.  Moreover, such twists are naturally investigated in the
context of the squares, not the triangles, of the octahedral vertex
figure we write Assessors on:  the three orthogonal Catamarans,
then, instead of the four touching-only-at-the-vertices Sails.
That's because opposite sides of a Catamaran twist to Assessors in
the same box-kite, so that each Catamaran lets one twist to two
different box-kites -- with the terminal Catamaran, in each case,
being further twistable into the box-kite you \textit{didn't} twist
to in the first instance.  As shown in the ``Twisted Sister'' and
``Royal Hunt'' diagrams of [4], these triple transforms can be
represented in their own Fano planes, with the indices placed on
their loci now corresponding to the \textit{strut constants} of a
septet of box-kites.

Each Catamaran comprises the pathways connecting 4 Assessors --
meaning it doesn't connect up with either term of the third strut in
its box-kite.  It is not hard to see that the strut constant of the
box-kite one twists to is equal to the strut-opposite of the term
which completes the L-trip of the edge being twisted in the first
place.  Hence, any L-index term on a Sedenion box-kite corresponds
to the strut constant of another such box-kite one can twist to.
This suggests expanding the meaning of ``twist product'' to embrace
pairings which share a strut rather than an edge.  For, if we allow
this, we can then treat the third strut orthogonal to the square
hosting twists as the ``mast'' of the Catamaran, giving us an
expanded sense of this latter term which allows us a major
simplification:  instead of thinking of the Sedenions' ZDs as
distributed among 7 distinct box-kites, we can see them all included
in one ``embroidered'' box-kite diagram, which we call a
\textit{brocade}.  Each of the 12 box-kite edges allows twists to a
pair of different Assessors -- let's say (A, b) and (B, a), in the
box-kite with ${\bf S} = c_{opp}$ = d.  More, the $({\bf S, X})$
pair -- which we can think of as \textit{in the box-kite's center}
-- can be ``twisted'' with all 6 Assessors in the original box-kite
to yield 12 more.  We therefore have 6 + 24 + 12 = the total set of
``42 Assessors'' in the Sedenions, all representable, on any one of
the 7 component box-kites, as a unitary ``brocade.''

It would be nice to be able to generalize the ``brocade'' notion so
as to reduce the number of basic structures in higher-order
contexts:  in the Pathions, for instance, there are 77 box-kites,
all but 21 of which are ``Type I,'' with 21 of those coming in
sand-mandala triples, 7 forming the ${\bf S} = 8$ ``Atlas," plus 7
ZP's and $3 \cdot 7$ ``strongboxes'' (so called, because these
low-${\bf S}$ box-kites contain ``pieces of 8'') completing the
collection. But if we also count in the $4 \cdot 7 = 28$ ``missing''
box-kites for high ${\bf S}$, we can collapse our head count from
105 box-kite-like structures to 15 brocades.  Miming the Sedenion
situation, the 7 ZP's form the simplest; the 7 Sand-Mandala trios
intermingle with the Atlas septet and the 21 strongboxes to make 7
more brocades; and the 21 Type II box-kites twist into each other
(to fill out one Catamaran in each) and into the ``hidden
box-kites'' linked with high ${\bf S}$ (filling out two more
Catamarans per Type II instance), yielding up the final set of 7
brocades.  (We note that the Type II situation is not as mysterious
as it might appear, once we recall the ``slipcover proof'' logic of
Part I, Section 5:  with 2 of 3 strut triplets being reversed,
``tugging'' on a Type II's Fano will tend to send a
\textit{reversed} arrow onto an edge 4 times out of 6 -- meaning
that, in all such cases, the corollary to Theorem 7, and hence the
theorem itself, will fail, thereby explaining the ``why'' of
``missing'' box-kites!)

We gain the generalized ``brocade'' simplification at a very small
price: relaxing the notion of ``twist product'' to embrace source
and target L- and U- index pairs which aren't necessarily
zero-divisors within the context of the ${\bf G}$ at hand. But this
is an investment which pays dividends, since it allows us to use
Type II structures as ``middlemen'' to facilitate studying the
``hidden box-kite'' substructures of the meta-fractal ``white
space'' in high-${\bf S}$ ET's.  Given the semiotic and semantic
importance of ``ZD-free'' structures (recall that our transcription
of Petitot's analysis of Greimas' ``Semiotic Square'' into
zero-divisor theory is based on ZD-free strut opposites), we can
expect a richness of results based on Catamaran study that should at
least equal that we are conducting based on Sails.  (For a ``coming
attraction,'' interested readers should see the online Powerpoint
slide-show linked with our NKS 2007 presentation [11], which will
play a role with respect to our forthcoming and similarly named
monograph, ``Voyage by Catamaran,'' akin to that our NKS 2006
slide-show did for the theorem/proof exposition you are currently
reading.)

\section*{References}
\begin{description}

\item \verb|[1]| Robert P. C. de Marrais, ``Placeholder
Substructures:  The Road from NKS to Small-World, Scale-Free
Networks Is Paved with Zero-Divisors,'' http:// \newline
wolframscience.com/conference/2006/
presentations/materials/demarrais.ppt  (Note:  the author's surname
is listed under ``M,'' not ``D.'')

\item \verb|[2]| Robert P. C. de Marrais, ``Placeholder
Substructures I:  The Road From NKS to Scale-Free Networks is Paved
with Zero Divisors,'' \textit{Complex Systems}, 17 (2007), 125-142;
arXiv:math.RA/0703745.

\item \verb|[3]| Robert P. C. de Marrais, ``The 42 Assessors and
the Box-Kites They Fly,'' arXiv:math.GM/0011260.

\item \verb|[4]| Robert P. C. de Marrais, ``Presto! Digitization,'' arXiv:math.RA/0603281

\item \verb|[5]| Stephen Wolfram, \textit{A New Kind of Science},
(Wolfram Media, Champaign IL, 2002).  Electronic version at
http://www.wolframscience.com/nksonline.

\item \verb|[6]| Robert P. C. de Marrais, ``Flying Higher Than A
Box-Kite,'' \newline arXiv:math.RA/0207003.

\item \verb|[7]| Robert P. C. de Marrais, ``The Marriage of Nothing
and All:  Zero-Divisor Box-Kites in a `TOE' Sky'', in Proceedings of
the $26^{\textrm{th}}$ International Colloquium on Group Theoretical
Methods in Physics, The Graduate Center of the City University of
New York, June 26-30, 2006, forthcoming from Springer--Verlag.

\item \verb|[8]| Robert P. C. de Marrais, ``The `Something From
Nothing' Insertion Point'',
http://www.wolframscience.com/conference/2004/presentations/materials/
\newline
rdemarrais.pdf

\item \verb|[9]| Robert P. C. de Marrais, ``Placeholder Substructures III:
A Bit-String-Driven `Recipe Theory' for Infinite-Dimensional
Zero-Divisor Spaces,'' \newline arXiv:0704.0112 [math.RA])

\item \verb|[10]| Benjamin Lee Whorf, \textit{Language, Thought, and
Reality}, edited by John B. Carroll (M.I.T. Press, Cambridge MA,
1956).

\item \verb|[11]| Robert P. C. de Marrais, ``Voyage by Catamaran:
Long-Distance Semantic Navigation, from Myth Logic to Semantic Web,
Can Be Effected by Infinite-Dimensional Zero-Divisor Ensembles,''
wolframscience.com/
\newline
conference/2007/presentations/materials/demarrais.ppt (Note:  the
author's surname is listed this time American style, under ``D,''
not ``M.'')

\end{description}

\end{document}